\documentclass[journal]{IEEEtran}

\usepackage{xcolor}
\usepackage{cite}
\usepackage{soul}

\ifCLASSINFOpdf
   \usepackage[pdftex]{graphicx}
   \graphicspath{{./Images/}}
   \DeclareGraphicsExtensions{.pdf,.jpeg,.png}
\else
\fi

\usepackage[cmex10]{amsmath}
\usepackage{amsmath,amsthm,mathtools}
\usepackage{amssymb}
\ifCLASSOPTIONcompsoc
  \usepackage[caption=false,font=normalsize,labelfont=sf,textfont=sf]{subfig}
\else
  \usepackage[caption=false,font=footnotesize]{subfig}
\fi

\usepackage{siunitx}
\usepackage{ifthen}
\usepackage{calc}
\usepackage{makeidx}
\usepackage{nomencl}
\usepackage{multirow}
\usepackage{footnote}

\usepackage{euscript}

\usepackage{flushend}

\makenomenclature
\nomenclature[aa]{$\mathcal{B}$}{Set of buses.}
\nomenclature[ab]{$\mathcal{G}$}{Set of generators $g$.}
\nomenclature[ac]{$\mathcal{M}$}{Set of demand aggregators $m$.}
\nomenclature[ad]{$\mathcal{L}$}{Set of lines: $\mathcal{L} \subseteq \mathcal{B} \times \mathcal{B}$.}
\nomenclature[ae]{$\mathcal{R}$}{Set of regions $r$.}
\nomenclature[af]{$\mathcal{H}$}{Set of time slots $h$.}
\nomenclature[ag]{$\underline{x}$/$\overline{x}$}{Minimum/maximum value of $x$.}
\nomenclature[ah]{$r_g^{+/-}$}{Ramp-up/down rate of supplier $g$.}
\nomenclature[ai]{$\tau_g^{\text{u}/\text{d}}$}{Minimum up/down time of supplier $g$.}
\nomenclature[aj]{$B_{i,j}$}{Susceptance of line between buses $i$ and $j$.}
\nomenclature[ak]{$\theta_{i,h}$}{Voltage angle at bus $i$ in $h$.}
\nomenclature[al]{$p_{l,h}^{i,j}$}{Power flow on line $i$-$j$: $p_{l,h}^{i,j} = B_{i,j}(\theta_{i,h}-\theta_{j,h})$.}
\nomenclature[am]{$\Delta p_{l,h}^{i,j}$}{Power loss on line $i$-$j$.}
\nomenclature[an]{$p_{g,h}$}{Generated power of supplier $g$ in $h$.}
\nomenclature[ao]{$p_{\text{b},h}^m$}{Battery storage power of aggregator $m$ in $h$.}
\nomenclature[ap]{$e_{\text{b},h}^m$}{Battery storage state of charge of aggregator $m$ in $h$.}
\nomenclature[aq]{$p_{\text{d},h}^{\text{\text{inf},m}}$}{Inflexible demand of aggregator $m$ in $h$.}
\nomenclature[ar]{$p_{\text{d},h}^{\text{flx},m}$}{Flexible demand of aggregator $m$ in $h$.}
\nomenclature[as]{$p_{\text{d},h}^{\text{u},m}$}{Underlying prosumer demand of aggregator $m$ in $h$.}
\nomenclature[at]{$p_{\text{pv},h}^{m}$}{PV generation of aggregator $m$ in $h$.}
\nomenclature[au]{$p^{\text{res}}_{r,h}$}{Required reserve in region $r$ in the system in $h$.}
\nomenclature[av]{$s_{g,h}$}{Binary decision variable on/off status of supplier $g$.} 
\nomenclature[aw]{$u_{g,h}$}{Binary start-up decision variable of supplier $g$.} 
\nomenclature[ax]{$d_{g,h}$}{Binary shut-down decision variable of supplier $g$.} 
\nomenclature[ay]{$c_g^{\text{fix}/\text{var}}$}{Fix/variable cost of supplier $g$.}
\nomenclature[az]{$c_g^{\text{su}/\text{sd}}$}{Start-up/shut-down cost of supplier $g$.}
\nomenclature[ba]{$\eta_{\text{b}}$}{Battery round-trip efficiency.}

\nomenclature[bb]{$\lambda$}{Dual variable associated with an equality constraint.}
\nomenclature[bc]{$\mu$}{Dual variable associated with an inequality constraint.}
\nomenclature[bd]{$b$}{Binary variable introduced to maintain complementary slackness.}
\nomenclature[bd]{$M$}{Large positive number used in the MILP formulation of the lower-level KKT conditions.}

\makeatletter
\def\BState{\State\hskip-\ALG@thistlm}
\makeatother
\begin{document}

\title{Generic Demand Model Considering the Impact of Prosumers for Future Grid Scenario Analysis}

\author{Shariq~Riaz,~\IEEEmembership{Student Member,~IEEE,}
Hesamoddin~Marzooghi,~\IEEEmembership{Member,~IEEE,}
Gregor~Verbi\v{c},~\IEEEmembership{Senior Member,~IEEE,}
Archie~C.~Chapman,~\IEEEmembership{Member,~IEEE,}
and~David~J.~Hill,~\IEEEmembership{Life Fellow,~IEEE}			
\thanks{Shariq~Riaz, Gregor~Verbi\v{c}, Archie~C.~Chapman, and David~J.~Hill are with the School of Electrical and Information Engineering, The University of Sydney, Sydney, New South Wales, Australia. e-mails: {shariq.riaz, gregor.verbic, archie.chapman, david.hill}@sydney.edu.au.}
\thanks{Hesamoddin~Marzooghi is with School of Engineering and Technology, Central Queensland University (CQ University), Brisbane, Australia. e-mail: h.marzooghi@cq.edu.au.}
\thanks{Shariq~Riaz is also with the Department of Electrical Engineering, University of Engineering and Technology Lahore, Lahore, Pakistan.}
\thanks{David~J.~Hill is also with the Department of Electrical and Electronic Engineering, The University of Hong Kong, Hong Kong.} 
\thanks{}}

\maketitle

\begin{abstract}
The increasing uptake of residential PV-battery systems is bound to significantly change demand patterns of future power systems and, consequently, their dynamic performance. In this paper, we propose a generic demand model that captures the aggregated effect of a large population of price-responsive users equipped with small-scale PV-battery systems, called \emph{prosumers}, for market simulation in future grid scenario analysis. The model is formulated as a bi-level program in which the upper-level unit commitment problem minimizes the total generation cost, and the lower-level problem maximizes prosumers' aggregate self-consumption. Unlike in the existing bi-level optimization frameworks that focus on the interaction between the wholesale market and an aggregator, the coupling is through the prosumers' demand, not through the electricity price. That renders the proposed model market structure agnostic, making it suitable for future grid studies where the market structure is potentially unknown. As a case study, we perform steady-state voltage stability analysis of a simplified model of the Australian National Electricity Market with significant penetration of renewable generation. The simulation results show that a high prosumer penetration changes the demand profile in ways that significantly improve the system loadability, which confirms the suitability of the proposed model for future grid studies.
\end{abstract}


\IEEEpeerreviewmaketitle

\printnomenclature[1.0 cm]

\section{Introduction}

\IEEEPARstart{P}{ower} systems are undergoing a major transformation driven by the increasing uptake of variable renewable energy sources (RES). At the demand side, the emergence of cost-effective ``behind-the-meter'' distributed energy resources, including on-site generation, energy storage, electric vehicles, and flexible loads, and the advancement of sensor, computer, communication and energy management technologies are changing the way electricity consumers source and consume electric power. 
Indeed, recent studies suggest that rooftop PV-battery systems will reach retail price parity from 2020 in the USA grids and the Australian National Electricity Market (NEM) \cite{EPRI}. 
A recent forecast by Morgan Stanley has suggested that the uptake can be even faster, by boldly predicting that up to 2 million Australian households could install battery storage by 2020 \cite{MorganStanley}. This has been confirmed by the Energy Networks Australia and the Australian Commonwealth Scientific and Industrial Research Organisation (CSIRO) who have estimated the projected uptake of solar PV and battery storage in 2050 to be \SI{80}{\giga \watt} and \SI{100}{GWh} \cite{ENA_CSIRO2006}, which will represent between 30\%--50\% of total demand, a scenario called ``Rise of the Prosumer'' \cite{CSIRO}. Here, the \emph{prosumer} they refer to is a small-scale (residential, commercial and small industrial) electricity consumer with on-site generation. A similar trend has been observed in Europe as well \cite{CEDelft2016}.
Given this, it is expected that a large uptake of demand-side technologies will significantly change demand patterns in future grids\footnote{We interpret a \textit{future grid} to mean the study of national grid type structures with the above-mentioned transformational changes for the long-term out to 2050.}, which will in turn affect their dynamic performance. 

Existing future grid feasibility studies \cite{Energy2010, Elliston2013, Budischak2013, Mason2010, AEMORES} typically use conventional demand models, possibly using some heuristics to account for the effect of emerging demand-side technologies, and the synergies that may arise between them. 
They also assume specific market arrangements by which RES are integrated into grid operations. 
The challenge associated with future grid planning is that the grid structure and the regulatory framework, including the market structure, cannot be simply assumed from the details of an existing one. Instead, several possible evolution paths need to be accounted for. Future grid planning thus requires a major departure from conventional power system planning, where only a handful of the most critical scenarios is analyzed. To account for a wide range of possible future evolutions, \emph{scenario analysis} has been proposed in many industries, e.g. in finance and economics \cite{LearningFromTheFuture_1998}, and in energy \cite{Foster2013, eHighway2015, NREL2012}. As opposed to power system planning where the aim is to find an optimal transmission and/or generation expansion plan, the aim in scenario analysis is to analyze possible evolution pathways to inform policymaking. Given the uncertainty associated with long-term projections, the focus of future grid scenario analysis is limited to the analysis of what is technically possible, although it might also consider an explicit costing \cite{Elliston2016}.
Our work is part of the Future Grid Research Program funded by the CSIRO, whose aim is to explore possible future pathways for the evolution of the Australian grid out to 2050 by looking beyond simple balancing. To this end, a comprehensive modeling framework for future grid scenario analysis has been proposed in \cite{Marzooghi2014}, which includes a market model, power flow analysis, and stability analysis. The demand model, however, assumes that the users are price-takers, which does not properly capture the aggregated effect of prosumers on the demand profile, as discussed later.

\subsection{Related Work}

Due to the influence of a demand profile on power systems performance and stability, recent studies have attempted to integrate the aggregated impact of prosumers into the demand models \cite{Claudio, Bruninx2012, Zungo13, Mathieu2014, Megel2015, GonzalezVaya2015}. The focus, however, is usually on scheduling of particular emerging demand-side technologies, e.g. HVAC \cite{Claudio, Bruninx2012, Zungo13}, flexible loads \cite{Mathieu2014}, PV-battery systems \cite{Megel2015}, and plug-in electrical vehicles (PEVs) \cite{GonzalezVaya2015}. Most of these modeling approaches assume an existing market structure, with the impact of prosumers incorporated by allowing demand and supply to interact in some limited or predefined ways. Specifically, this is mainly done via three different approaches:
\begin{enumerate}
	\item{Only the supply-side is modeled physically while prosumers are considered by a simplified representation of demand-side technologies. In \cite{Mathieu2014}, flexible loads' effects on reserve markets are analyzed by modeling prosumers with a tank model; however, the reserve market is greatly simplified. In \cite{Claudio}, prosumers are represented by a price-elasticity matrix, which is used to model changes in the aggregate demand in response to a change in the electricity price, and are acquired from the analysis of historical data.}
	\item{Demand-side technologies are physically modeled while a simplified representation of supply-side is employed. For instance, in \cite{Megel2015}, the supply-side is represented by an electricity price profile.}
	\item{Both supply and demand sides can be modeled physically and optimized jointly, as in \cite{Bruninx2012, GonzalezVaya2015}, which can produce more realistic results. For example, the study in \cite {GonzalezVaya2015} integrates the aggregated charging management approaches for PEVs into the market clearing process, with a simplified representation of the latter.}
\end{enumerate}
Although the above models have shown their merits, they are dependent on specific practical details such as the electricity price or the implementation of a mechanism for demand response (DR) aggregation, which limits their usefulness for future grid scenario analysis where the detailed market structure is potentially unknown.

Against this backdrop, the paper proposes a principled method for generic demand modeling including the aggregated effect of prosumers. The model is formulated as a bi-level program in which the upper-level unit commitment problem minimizes the total generation cost, and the lower-level problem maximizes the \emph{aggregate} prosumers' self-consumption. 
In more detail, the lower-level objective is motivated by the emerging situation in Australia, where rooftop PV owners are increasingly discouraged from sending power back to the grid due to very low PV feed-in-tariffs versus increasing retail electricity prices. In this setting, an obvious cost-minimizing strategy is to install small-scale battery storage, to maximize self-consumption of local generated energy and offset energy used in peak pricing periods. Similar tariff settings appear likely to occur globally in the near future, as acknowledged in \cite{GP}.
Moreover, self-consumption within an aggregated block of prosumers is a good approximation of many likely behaviors and responses to other future incentives and market structures, such as (peak power-based) demand charges, capacity constrained connections, virtual net metering across connection points, transactive energy and local energy trading, and a (somewhat irrational) desire for self-reliance.

A key difference from existing bi-level optimization frameworks is that in our formulation, the levels are coupled through the prosumers' demand, not through the electricity price. In contrast, other models, which focus on the interaction between an aggregator and the prosumers \cite{Zungo13} or the aggregator and the wholesale market \cite{GonzalezVaya2015}, couple the levels through prices. These approaches essentially define a market structure, that is, a pricing rule to support an outcome. In contradistinction, our proposed model is market structure agnostic. That is, it implicitly assumes that an efficient mechanism for demand response aggregation is adopted, with prices determined by that unspecified mechanism, which support the outcomes computed by our optimization framework.

The paper builds on our previous work \cite{SEGAN2015, Riaz2016}. In \cite{Riaz2016}, we have assumed that a prosumer aggregation represents a homogeneous group of loads; that is, we have assumed that they all behave in the same way and have the same capacity, and are allowed to send power back to the grid. In the absence of an explicit transmission pricing, this can create perverse outcomes, such as power exchange between aggregators located in different parts of the network. In the model proposed in this paper, the aggregators are not allowed to send the power back to the grid, which better models the assumption of self-consumption within an aggregated block of prosumers.
The model proposed in \cite{SEGAN2015} is similar to the model in this paper, however, the prosumer battery storage is modeled implicitly, which requires a heuristic search to capture the prosumer behavior.

The remainder of the paper is organized as follows: Section II presents the proposed modeling framework. In Section III, the efficacy of the proposed framework is demonstrated on a simplified 14-generator network model of the NEM with a significant RES penetration. Finally, Section IV concludes.

\section{Generic Demand Model Considering the Impact of Prosumers}
Generic demand models are essential for power system studies. They are commonly used to reflect the aggregated effect of numerous physical loads \cite{Kundur, VanCutsem}. Conventional demand models only account for the accumulated effect of independent load changes and some relatively minor control actions, that is, they are by necessity simplified to not include all details of the loads. In the following, we explain the motivation behind this work and the modeling assumptions.

\subsection{Research Motivation and Modeling Assumptions}
\begin{figure} 
\centering
\includegraphics[width=8.0cm, keepaspectratio]{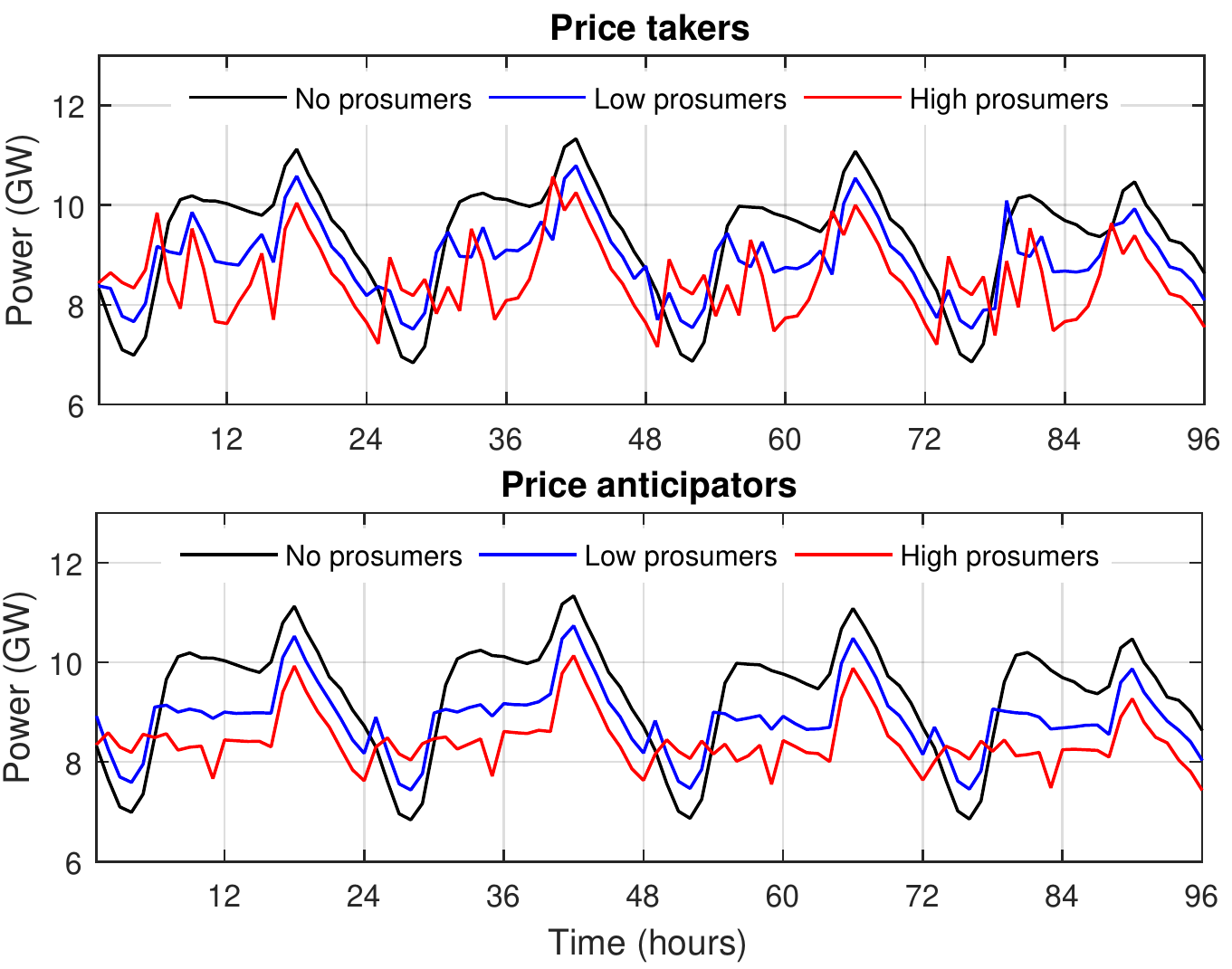}
\caption{Operational demand with different penetrations of price-taking (top) and price-anticipating prosumers (bottom).}
\label{figure:The effect of different DR scenarios on the demand profile}
\end{figure}

The main purpose of developing generic demand models is to provide accurate dispatch decisions for balancing and stability analysis of future grid scenarios. 
Given the uncertainty associated with future grid studies, the modeling framework should be market structure agnostic, and capable of easy integration of various types and penetrations of emerging demand-side technologies. To this effect, we make the following assumptions:
\begin{enumerate}
	\item{The loads are modeled as price anticipators. 
	It is well understood that price-taking load behavior, which simply responds to a given price profile, can result in load synchronization, i.e. all users move their consumption to a low-price period, resulting in an inefficient market outcome. 
	In contrast, price anticipatory loads influence the electricity price by playing ``a game'' with the wholesale market. The game is captured by the proposed bi-level model. 
	Specifically, the price-anticipating assumption implies that the aggregate effect of the prosumers is to change the market clearing prices and quantities, and, moreover, that the prosumers have a model of this effect. 
	Given this, the prosumer aggregation bids follow an equilibrium strategy, with an accurate expectation about the response of the market; this is the standard reasoning behind a Cournot or Stackelberg game formulation corresponding to a bi-level optimization problem. This equilibrium strategy can be thought of as being generated by the following iterative process: First, the market operator creates a price profile by clearing the market based on the predicted demand. Second, the prosumers (and other price-anticipatory participants) respond by shifting their consumption to cheaper time slots. This gives a new demand profile, and the market operator clears the market again, and the process repeats until convergence. 
	Note that the proposed model does not specify an iterative mechanism, but rather, it encodes the optimality conditions for any price profile for the prosumers in Karush-Kuhn-Tucker (KKT) conditions, 
	and in this way it captures the outcome resulting from the price-anticipatory prosumer behavior. 
	In our previous work \cite{Hesam2014, Marzooghi2014}, we modeled the loads as price-takers, inspired by the smart home concept \cite{Tischer2011}, in which the loads respond to the electricity price to minimize energy expenditure. The study has shown that with large penetration of price-taking prosumers, the marginal benefit might become negative when secondary peaks are created due to the load synchronization (see Fig.~\ref{figure:The effect of different DR scenarios on the demand profile} showing the operational demand with different penetrations of prosumers). Therefore, demand response aggregators (henceforth simply called aggregators) have started to emerge to fully exploit the demand-side flexibility. To that effect, the model implicitly assumes an efficient mechanism for demand response aggregation (an interested reader might refer to \cite{AGD} for a discussion on practical implementation issues). However, specific implementation details, like price structure or the division of the profit earning by the aggregated collection of prosumers, are not of explicit interest in the proposed model.}
	\item{The demand model representing an aggregator consists of a large population of prosumers connected to an unconstrained distribution network who collectively maximize self-consumption (made possible by an efficient internal trading and balancing mechanism).}
	\item{Aggregators do not alter the underlying power consumption of the prosumers. That is, except for battery losses, the total power consumption before and after aggregation remains the same; however, the grid power intake profile does change, as a result of the prosumers using batteries to maximize self-consumption.}
	\item{Prosumers have smart meters equipped with home energy management (HEM) systems for scheduling of the PV-battery systems. Also, a communication infrastructure is assumed that allows a two-way communication between the grid, the aggregator and the prosumers, facilitating energy trading between prosumers in the aggregation. }
\end{enumerate}
These assumptions appear to be appropriate for scenarios arising in time frame of several decades into the future. 

\subsection{Bi-level Optimization Framework}
In the model, we are specifically interested in the aggregated effect of a large prosumer population on the demand profile, assuming that the prosumers collectively maximize their self-consumption. Given that the objective of the wholesale market is to minimize the generation cost, the problem exhibits a bi-level structure. In game theory, such hierarchical optimization problems are known as Stackelberg games. They can be formulated as bi-level mathematical programs of the form \cite{Zungo13}:
\begin{align*}
	\mathop{\operatorname{minimize}}\limits_{ \textbf{x}, \textbf{y} } \quad & \Phi \left( \textbf{x}, \textbf{y} \right) \\
	\operatorname{subject\,to} \quad & \left( \textbf{x}, \textbf{y} \right) \in \mathcal{Z} \\	
	& \textbf{y} \in \mathcal{S} = \mathop{\operatorname {arg\,min}}\limits_{\textbf{y}} \{ \Omega \left( \textbf{x}, \textbf{y} \right): \textbf{y} \in \mathcal{C} \left( \textbf{x} \right) \} 
\end{align*}
where $\textbf{x} \in \mathbb{R}^n$, $\textbf{y} \in \mathbb{R}^m$, are decisions vectors, and $\Phi \left( \textbf{x}, \textbf{y} \right): \mathbb{R}^{n+m} \rightarrow \mathbb{R}$ and $\Omega \left( \textbf{x}, \textbf{y} \right): \mathbb{R}^{n+m} \rightarrow \mathbb{R}$ are the objective functions of the upper- and the lower-level problems, respectively. $\mathcal{Z}$ is the joint feasible region of the upper-level problem and $\mathcal{C(\textbf{x})}$ the feasible region of the lower-level problem induced by $\textbf{x}$.
In the existing market models that adopt a hierarchical approach, the coupling variable $\textbf{y}$ is
the electricity price (e.g. \cite{GonzalezVaya2015, Zungo13}). That is, the upper-level (the wholesale market in our case) determines the price schedule, while the lower-level (the aggregator acting on behalf of the prosumers), optimizes its consumption based on this price schedule. 

\begin{figure} 
\centering
\includegraphics[width=7.0cm, keepaspectratio]{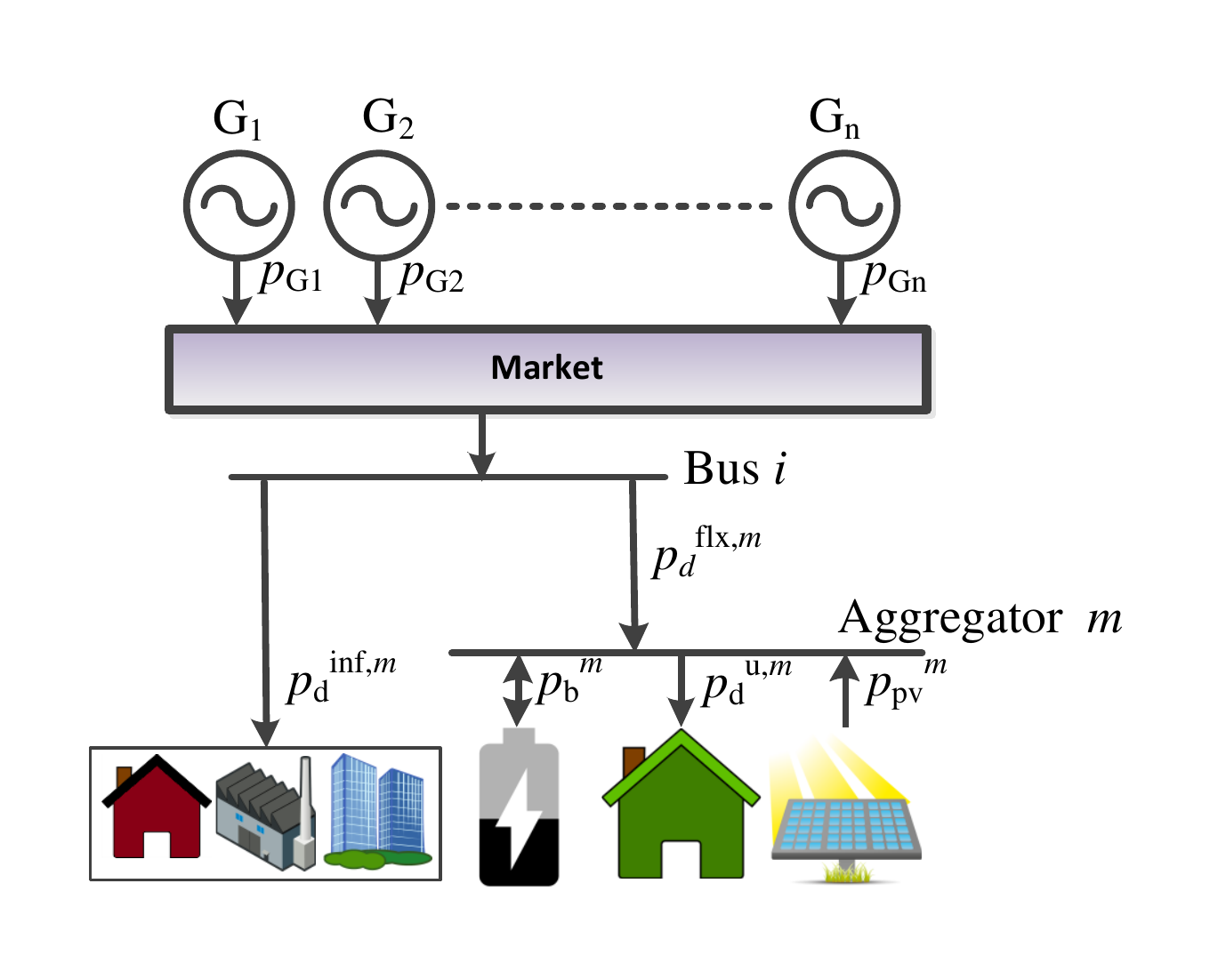}
\caption{Structure of the proposed modeling framework.}
\label{figure:LoadModel}
\end{figure}

Fig. \ref{figure:LoadModel} shows the structure of the proposed modeling framework. The demand model consists of two parts: (i) \emph{inflexible} demand, $p_{\text{d}}^{\text{inf},m}$, with a fixed demand profile, representing large industrial loads and loads without flexible resources; and (ii) \emph{flexible} demand, $p_{\text{d}}^{\text{flx},m}$, comprising a large population of prosumers who collectively maximize self-consumption. Note that not every bus in the system has a load connected to it, hence the distinction between an aggregator $m \in \mathcal{M}$ and bus $i \in \mathcal{B}$.
Unlike in most existing studies, the interaction between the wholesale market and the aggregators in our model is through the demand profile of the aggregator, $p_{\text{d}}^{\text{flx},m}$.
Note that in contradistinction to the price-taking assumption when the electricity price is known in advance, now the collective action of the prosumers affects the wholesale market dispatch, which is the salient feature of the proposed model.

\subsection{Upper-level Problem (Wholesale Market)} \label{ULProblem} 
To emulate the market outcome, the upper-level problem is cast as a unit commitment (UC) problem aiming to minimize the generation cost:
\begin{align}\label{eq:UC_utility}
	\mathop{\operatorname{minimize}} \limits_{s,u,d,p,\theta} \sum_{h \in \mathcal{H}} \! {\sum_{g \in \mathcal{G}}} \Big( 
	c_g^\text{fix} s_{g,h} \! + \! c_g^\text{su} u_{g,h} \! + \! c_g^\text{sd} d_{g,h} \! + \! c_g^\text{var} p_{g,h} \Big),
\end{align}
where $s_{g,h},u_{g,h},d_{g,h} \in \{0,1\}$, $p_{g,h}\in\mathbb{R}_{+}$, $\theta_{i,h}\in\mathbb{R}$ are the decision variables of the problem. The problem is subject to the following constraints: 
\begin{align}\label{LC2}
	& \sum_{g \in \mathcal{G}_i}{p_{g,h}} = \! \sum_{m \in \mathcal{M}_i} \! (p_{\text{d},h}^{\text{inf},m}+p_{\text{d},h}^{\text{flx},m})+\sum_{l \in \mathcal{L}_i}(p_{l,h}+\Delta p_{l,h}), \\
	& |B_{i,j}(\theta_{i,h}-\theta_{j,h})| \leq \overline{p}_{l},\label{LC4} \\
	& \underline{p}_{g} s_{g,h} \leq p_{g,h} \leq \overline{p}_{g} s_{g,h}, \label{LC8} \\
	& u_{g,h}-d_{g,h}=s_{g,h}-s_{g,h-1},\label{LC9} \\
	& \sum\nolimits_{g_{\text{synch}}\in \mathcal{R}}({\overline{p}_{g}} s_{g,h} - p_{g,h}) \geq p^{\text{res}}_{r,h},\label{LC10} \\
	& u_{g,h} + \sum\nolimits_{\tilde{h}=0}^{\tau_g^{\text{u}}-1} d_{g,h+\tilde{h}} \leq 1,	\label{LC11} \\
	& d_{g,h} + \sum\nolimits_{\tilde{h}=0}^{\tau_g^{\text{d}}-1} u_{g,h+\tilde{h}} \leq 1,	\label{LC12} \\
	& -r_g^- \leq p_{g,h} - p_{g,h-1} \leq r_g^+,	\label{LC13}\\
	&\mathop{\operatorname {arg\,min}} \limits_{p_{\text{d}}^{\text{flx}}, p_{\text{b}}, e_{\text{b}} } \left\lbrace  \sum_{h \in \mathcal{H}} p_{\text{d},h}^{\text{flx},m} \ \operatorname{subject\,to} \  \eqref{LC15}\textrm{-}\eqref{LC18} \right\rbrace , \label{LC20} 
\end{align}
where \eqref{LC2} is the power balance equation at each bus $i$ in the system\footnote{Note that the flexible demand of each aggregator $m$, $p_{\text{d},h}^{\text{flx},m}$, couples the upper-level (wholesale market) problem with each of the $m$ lower-level (aggregator) problems.}, with $\mathcal{G}_i$, $\mathcal{M}_i$, $\mathcal{L}_i$ representing respectively the sets of generators, aggregators and lines connected to bus $i$, and $p_{\text{d},h}^{\text{inf},m}$, $p_{\text{d},h}^{\text{flx},m}$, $p_{l,h}^{i,j}$ and $\Delta p_{l,h}^{i,j}$ representing respectively the inflexible and flexible demand of aggregator $m$, line power and line power loss (assumed to be 10\% of the line flow) on each line connected to bus $i$;
\eqref{LC4} represents line power limits;
\eqref{LC8} limits the dispatch level of a generating unit between its respective minimum and maximum limits;
\eqref{LC9} links the status of a generator unit to the up and down binary decision variables;
\eqref{LC10} ensures sufficient spinning reserves are available in reach region of the grid;
\eqref{LC11} and~\eqref{LC12} ensure minimum up and minimum down times of the generators; 
\eqref{LC13} are the generator ramping constraints; and
\eqref{LC20} is the constraint resulting from the prosumer aggregation optimization problem, as explained in Section~\ref{LLP}.

\subsection{Lower-level Problem (Aggregators)} \label{LLP}
The prosumer aggregation is formulated in the lower-level problem. 
The loads within an aggregator's domain are assumed homogeneous, 
which allows us to represent the total aggregator's demand with a single load model. 
The flexibility provided by the battery is only used to maximize self-consumption, with grid supply readily available.
This implies that the end-users' power consumption pattern is left unaltered, so that their comfort is not jeopardized. 

The coupling between the upper-level (wholesale) and the lower-level (retail) problem in the proposed model is through the power demand.
This removes the need to define the market in terms of a pricing mechanism or rule, 
and it follows that the electricity price is not explicitly shown in the optimization problem. This is why it is called a \emph{generic model}.

This approach stands in contrast to other existing bi-level formulations \cite{Zungo13},\cite{GonzalezVaya2015}, in which the loads and wholesale market are coupled through the electricity price. 
The prices generated by any market depend inherently on the specific pricing mechanism adopted. 
However, in practice, several different pricing rules can implement a desired outcome, 
including the widely-used uniform price reverse auction or nodal pricing mechanisms. 
For example, the electricity price could comprise the dual variables associated with the power balance constraint \eqref{LC2} and power flow constraints \eqref{LC4} of the upper-level problem, plus retail/aggregator and network charges. 
Moreover, a range of different pricing rules result in different retail/aggregator 
and network charges, with all supporting an efficient outcome.
However, by instead coupling the upper- and lower-level problems via power demand, 
our proposed model avoids the need to specify a particular pricing rule, 
which makes it market-structure agnostic (see also Assumptions 1 and 2 in Section II-A). 
Nonetheless, it is fair to assume that in any practical system, 
users will have access to a price forecast\footnote{Note that in a HEM problem \cite{Tischer2011}, a HEM system is an agent acting on behalf of a prosumer, and the electricity price is known ahead of time, resulting in a price-taking behavior.} and other bidders' historical behavior--but these are all market design-specific. 
In the absence of pricing rule details, 
the proposed optimization model does not require such information.

Specifically, the lower-level problem is formulated as follows for each demand aggregator $m \in \mathcal{M}$:
\begin{align}\label{LC14}
	& \mathop{\operatorname{minimize}}\limits_{ p_{\text{d}}^{\text{flx}}, p_{\text{b}}, e_{\text{b}} } \sum_{h \in \mathcal{H}} p_{\text{d},h}^{\text{flx},m},
\end{align}
where the decision variables of the problem are flexible demand $p_{\text{d}}^{\text{flx}}$, battery power $p_{\text{b}}$, and battery capacity $e_{\text{b}}$.
The problem is subject to the following constraints:
\begin{align}
	& p_{\text{d},h}^{\text{flx},m} = p_{\text{d},h}^{\text{u},m}-p_{\text{pv},h}^{m}+p_{\text{b},h}^m,\label{LC15}\\
	& e_{\text{b},h}^m = \eta_{\text{b}}^m e_{\text{b},h-1}^m + p_{\text{b},h}^m, \label{LC19} \\
	& \underline{p}_{\text{b}}^m \leq p_{\text{b},h}^m \leq \overline{p}_{\text{b}}^m,\label{LC17}\\
	& \underline{e}_{\text{b}}^{m} \leq{e}_{\text{b},h}^{m} \leq \overline{e}_{\text{b}}^{m},\label{LC18}
\end{align}
where \eqref{LC15} is the power balance equation; and 
\eqref{LC19}-\eqref{LC18} are the battery storage constraints.
Power $p_{\text{d},h}^{\text{u},m}$ is the underlying power demand of the prosumers. Note that according to the Assumption 3 in Section II-B, except for battery losses, the underlying power demand does not change, however the power intake from grid can. Finally, the KKT optimality conditions of the lower-level problem are added as the constraints to the upper-level problem, which reduces the problem to a single mixed integer linear program that can be solved using of-the-shelf solvers. Note that because the two levels interacts through a power, not through a price, unlike in \cite{Zungo13}, no linearization is required.

\section{Case Studies}
To showcase the efficacy of the model, we analyze steady-state voltage stability of a simplified model of the NEM with scenarios reflecting different prosumer penetrations.

\subsection{Model of the Australian National Electricity Market (NEM)}
\begin{figure} 
\centering
\includegraphics[width=7.9cm, keepaspectratio]{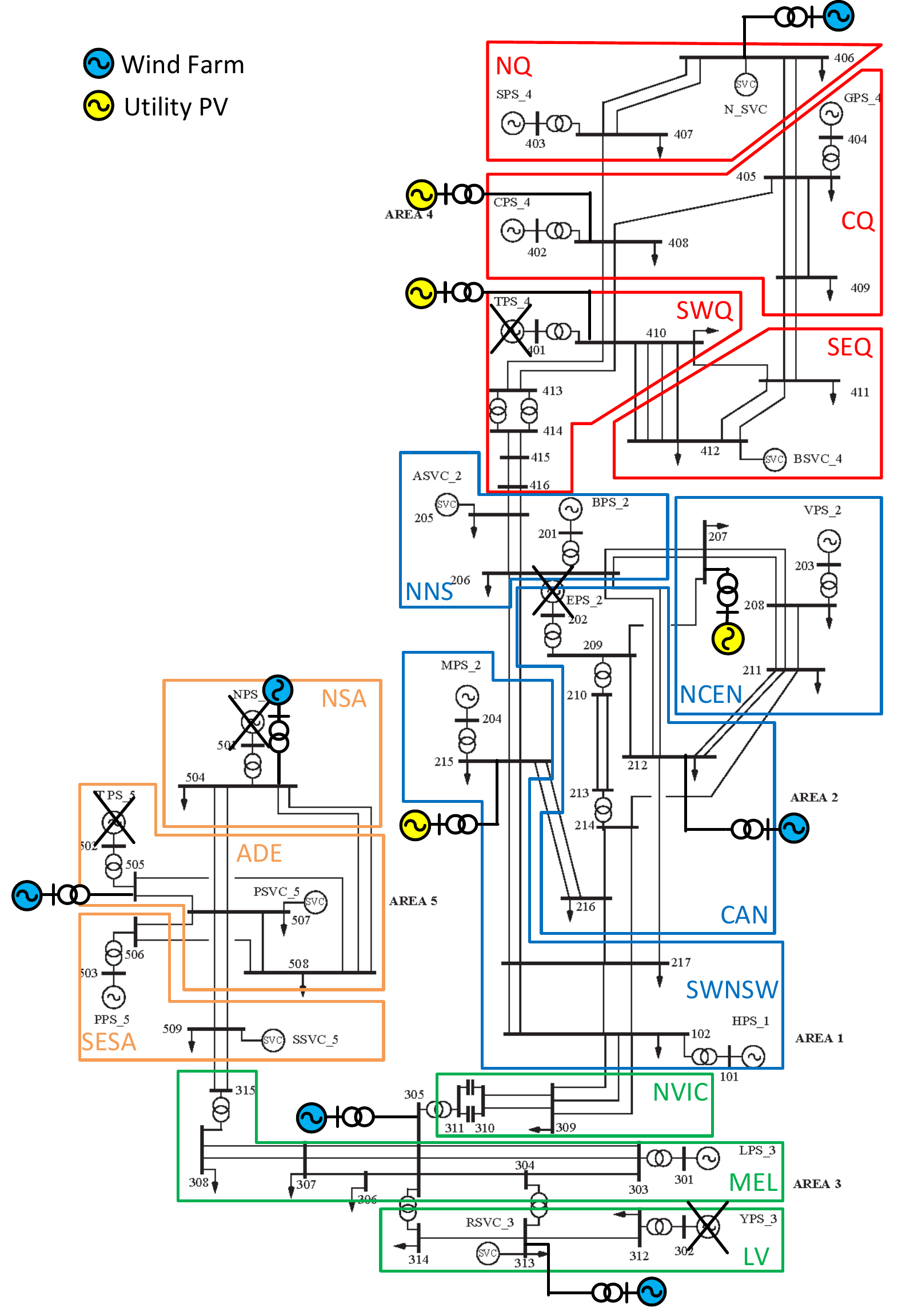}
\caption{Single-line diagram of the 14-generator model of the NEM with the 16 zones that define wind and solar traces \cite{AEMO2012A}.}
\label{fig:14-generator model of the NEM}
\end{figure}

The 14-generator IEEE test system shown in Fig. \ref{fig:14-generator model of the NEM} was initially proposed in \cite{Gibbard2010} as a test bed for small-signal analysis. The system is loosely based on the NEM, the interconnection on the Australian eastern seaboard. The network is stringy, with large transmission distances and loads concentrated in a few load centers. It consists of 59 buses, 28 loads, and 14 generators.
The test system consists of four areas representing the states of Queensland, New South Wales, Victoria and South Australia, and 28 aggregators, one for each load bus. 
The generator technologies and modeling assumptions follow \cite{SEGAN2015}. 
We consider two RES penetration rates. In the business as usual (BAU) scenario, the generation portfolio includes \SI{39.36}{\giga\watt} coal, \SI{5.22}{\giga\watt} gas, and \SI{2.33}{\giga\watt} hydro. In the high-RES scenario, \SI{40}{\percent} of the total demand is covered by variable RES. Inspired by two recent Australian 100\% renewables studies \cite{Energy2010, AEMORES}, part of coal generation is replaced with wind and utility PV using wind and solar traces from the AEMO's planning document \cite{AEMO2012A}, which results in \SI{28.94}{\giga\watt} coal, \SI{5.22}{\giga\watt} gas, \SI{2.33}{\giga\watt} hydro, \SI{21}{\giga\watt} wind, and \SI{12}{\giga\watt} utility PV. Given the deterministic nature of the model, we assume \SI{10}{\percent} reserves for each region in the system to cater for demand and RES forecast errors.
In market simulations, generators are assumed to bid according to their short-run marginal costs, while RESs bid at zero cost. 
Simulations are performed using a rolling horizon approach with hourly resolution assuming a perfect foresight. The optimization horizon is three days with a two-day overlap. Last, wind and solar generators are assumed to operate in a voltage control mode.

\subsection{Prosumer Scenarios}
We assume four different prosumer penetrations: zero, low, medium and high. With no prosumer penetration, the demand is assumed inflexible. For the other three scenarios, we assume that part of the demand is equipped with small-scale (residential and small commercial) PV-battery systems. The uptake of PV loosely follows a recent AEMO study \cite{AEMOTECH2015}. The PV capacities are respectively \SI{5}{\giga\watt}, \SI{10}{\giga\watt}, and \SI{20}{\giga\watt} for the low, medium and high uptake of prosumers (the existing penetration in the NEM is \SI{5}{\giga\watt}).
We consider three different amounts of storage: zero, \SI{2}{\kilo\watt\hour}, and \SI{4}{\kilo\watt\hour} of storage for \SI{1}{\kilo\watt} of rooftop PV.\footnote{A typical ratio in the NEM today is \SI{2}{h} of storage \cite{AEMOTECH2015}, however, in the future, this will likely increase due to the anticipated cost reduction.}
Hourly demand and PV traces are from the AEMO's planning document \cite{AEMO2012A}.

\subsection{Dispatch Results}
\begin{figure*}
 \center
 \includegraphics[width=\textwidth]{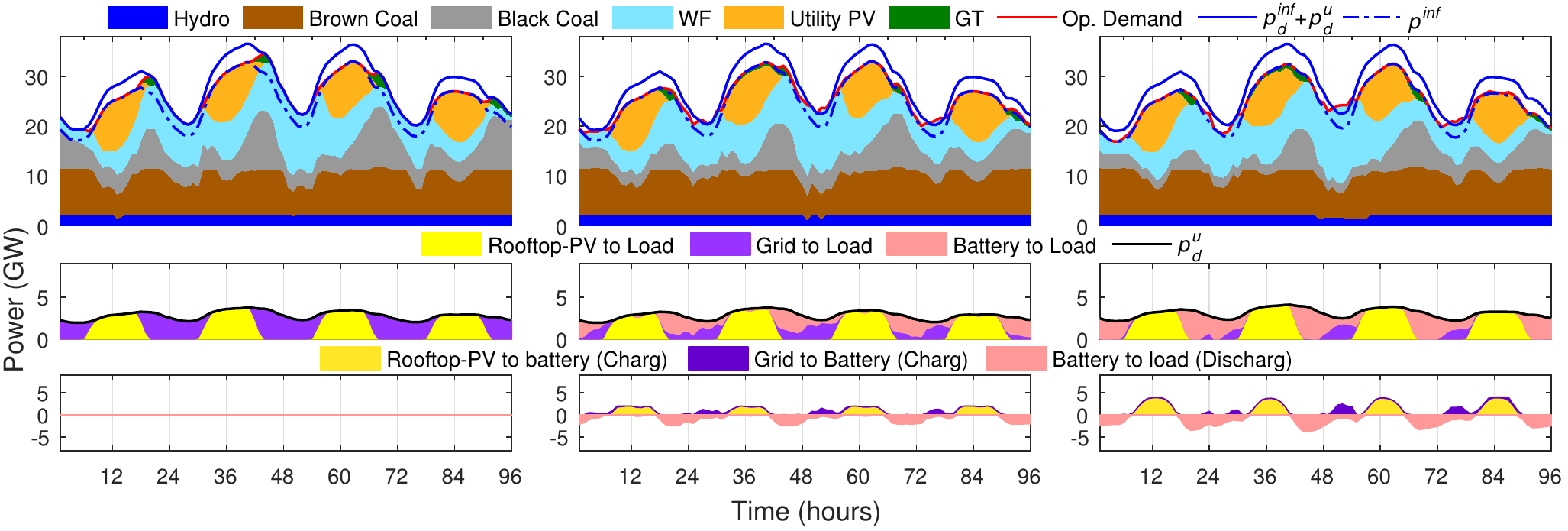}
 \caption{Dispatch results for a typical summer week with high demand (12-15 January) for a medium prosumer penetration with different amounts of storage: zero (left), \SI{2}{h} (middle) and \SI{4}{h} (right).}
 \label{fig:Gen_Mix_MP_0h2h4h}
\end{figure*}

\begin{figure*}
 \center
 \includegraphics[width=\textwidth]{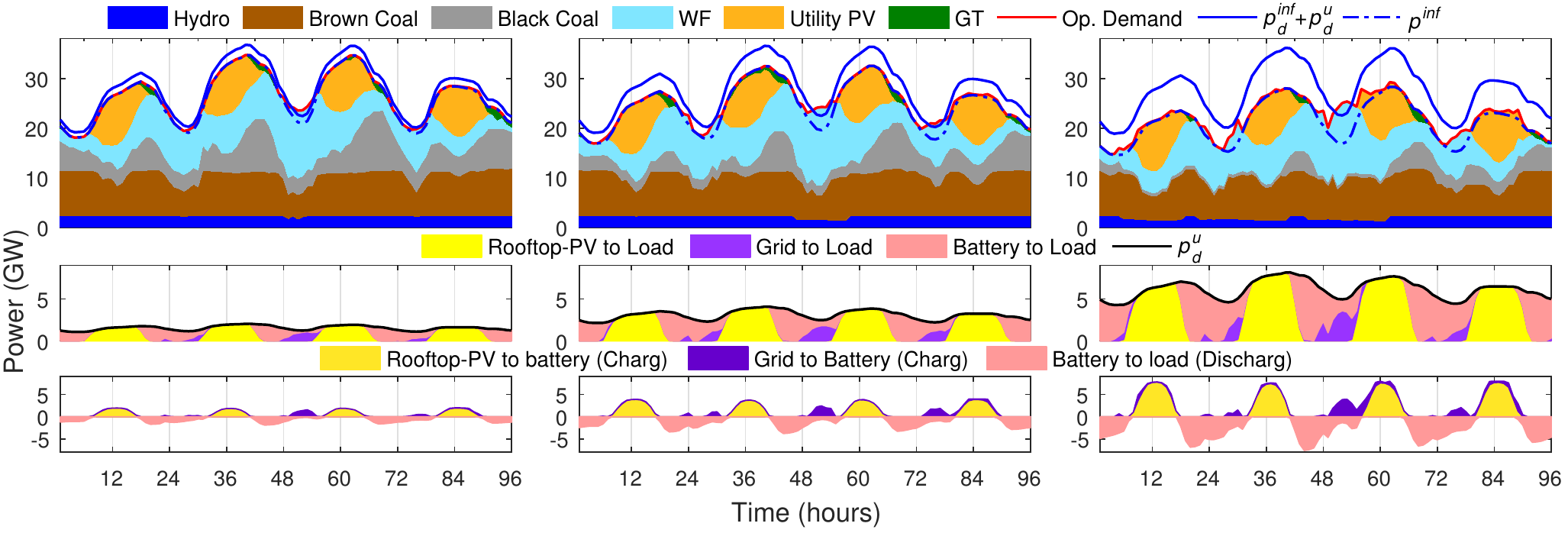}
 \caption{Dispatch results for a typical summer week with high demand (12-15 January) for different penetrations of prosumers with \SI{4}{h} of storage: zero (left), medium (middle), high (right).}
 \label{fig:Gen_Mix_LPMPHP_4h}
\end{figure*}

Dispatch results for a typical summer week with high demand (12-15 January) for a few representative scenarios are shown in Figs.~\ref{fig:Gen_Mix_MP_0h2h4h} and \ref{fig:Gen_Mix_LPMPHP_4h}. 
The figures show, respectively, generation dispatch results (top row), combined flexible demand of all aggregators (middle row), and a combined battery charging profile of all aggregators (bottom row).
Fig.~\ref{fig:Gen_Mix_MP_0h2h4h} shows results for a medium prosumer penetration with, respectively, zero, \SI{2}{h} and \SI{4}{h} hours of storage.
Fig.~\ref{fig:Gen_Mix_LPMPHP_4h} shows results for different prosumer penetrations (zero, medium, high) with \SI{4}{h} of storage. 
Observe that in all six cases peak demand occurs at mid-day due to a high air-conditioning load. After the sunset, however, the demand is still high, so gas generation is needed to cover the gap. In the available generation mix, gas has the highest short-run marginal cost, which increases the electricity price in late afternoon/early evening. 
The balancing results over the simulated year have revealed that the increased RES penetration in the renewable scenarios requires more energy from gas generation compared to the BAU scenario. This is due to RES intermittency, and the ramp limits of conventional coal-fired generation. An increased penetration of prosumers with higher amounts of storage, however, reduces the usage of gas due to a flatter demand profile.

Observe in Fig.~\ref{fig:Gen_Mix_MP_0h2h4h} how an increasing amount of storage increases prosumers' self-sufficiency. Without storage, the load is supplied by PV during the day, and the rest is supplied from the grid. When storage is added to the system, batteries are charged when electricity is cheap (mostly from rooftop PV during the day and from wind during the night) and discharged in late afternoon to offset the demand when the electricity is most expensive. Note that the plots in the bottom two rows show a combined load profile of all aggregators in the system, which explains why storage is seemingly charged and discharged simultaneously. 
Observe how high amounts of storage (rightmost columns in Figs.~\ref{fig:Gen_Mix_MP_0h2h4h} and \ref{fig:Gen_Mix_LPMPHP_4h}) flatten the demand profile. During the day, the flexible demand is supplied by rooftop PV, which \emph{reduces} the operational demand, while during the night, with sufficient wind generation, batteries are charged, which \emph{increases} the operational demand. This has a significant beneficial effect on loadability and voltage stability, as discussed in the next section.

\subsection{Loadability and Voltage Stability Results}
Dispatch results from the market simulations are used to perform a load flow analysis, which is then used in the assessment of loadability and voltage stability. In the analysis, only scenarios with \SI{4}{h} of storage were considered. The prosumer scenarios are thus called, according to the respective penetration rates, zero (ZP), low (LP), medium (MP), and high (HP).
Note that the market model only considers a simplified DC power flow with the maximum angle limit set to $30^{\circ}$. This can sometimes result in a non-convergent AC power flow in scenarios with a high RES penetration. The number of non-convergent hours is, respectively, 175, 37, 12, and 0, in scenarios ZP, LP, MP, and HP. An increased penetration of prosumers thus improves voltage stability, as explained in more detail later.

In loadability assessment, N-1 security is considered, so a contingency screening is performed first. We screened all credible N-1 contingencies to identify the most severe ones based on the maximum power transfer level \cite{Vaahedi1999}. Twenty most critical contingencies were selected for each hour of the simulated year. 

\subsubsection{Loadability calculation} 
To calculate system loadability (LDB), the load in the system is progressively increased using the results of the market dispatch as the base case. After each increase, load flow is calculated twice; first for the case after the last load increase without a contingency, and then again for each of the preselected critical contingencies. When the load flow does not converge for a particular contingency, the last convergent load flow solution without a contingency is considered the loadability margin.
We considered two different load increase patterns, where load and generation are increased uniformly, in proportion to the base case: (i) \textbf{NEM}: only load and generation in the NEM are increased; (ii) \textbf{SA/VIC}: only load in VIC and generation in SA are increased. The results are summarized in Table \ref{tab:loadability results}.
Comparing the BAU scenario and the renewable scenario with conventional demand (ZP), it can be seen that with the increased RES penetration, the average loadability margin over the simulated year in the demand increase scenario NEM is decreased from \SI{7.8}{\giga\watt} to \SI{2.5}{\giga\watt}. Similarly, the average loadability margin in the demand increase scenario SA/VIC is reduced from \SI{2.2}{\giga\watt} to \SI{0.8}{\giga\watt}.
With a high RES penetration, conventional synchronous generation is replaced by inverter-based generation with inferior reactive power support capability\footnote{For synchronous generation, a \SI{0.8}{} power factor is assumed. For RES, we used the reactive power capability curve for the generic GE Type IV wind farm model \cite{Clark2013}, in which the reactive power generation is significantly constrained close to the nominal active power generation. In practice, detailed planning studies are required before connection is granted, which can result in an increased requirement for reactive power support.}, which results in a reduced reactive power margin in the system and hence lower stability margin.
With an increased penetration of prosumers, the system loadability improves. Observe that the average system loadability margin in both load increase scenarios, NEM and SA/VIC, is increased from \SI{2.5}{\giga\watt} and \SI{0.8}{\giga\watt} for the renewable scenario with no prosumers (ZP) to \SI{7.1}{\giga\watt} and \SI{1.7}{\giga\watt} for high penetration of prosumers (HP), respectively, which indicates a considerable improvement in the system loadability margin. 
This is explained by a demand reduction when prosumer demand is supplied by rooftop PV. In the night hours, however, even with high prosumer penetration, the loadibility can be reduced when prosumers charge their batteries, as observed in Figs.~\ref{fig:Gen_Mix_MP_0h2h4h} and \ref{fig:Gen_Mix_LPMPHP_4h}.
The situation is further illustrated in Fig. \ref{figure:LDBSA} that compares the loadability margin for the load increase scenario NEM and the results of the modal analysis, discussed next.

\begin{table}
\centering
\caption{Loadability results}
\label{tab:loadability results}
\begin{tabular}{|c|c|c|c|c|c|c|}\hline
&\begin{tabular}[c]{@{}c@{}}LDB\\Cases\end{tabular}&\multicolumn{5}{c|}{Scenarios}\\\hline
&&BAU&ZP&LP&MP&HP\\\hline
\multirow{2}{*}{\begin{tabular}[c]{@{}c@{}}Avg. LDB\\margin (GW)\end{tabular}}&NEM&\begin{tabular}[c]{@{}c@{}}7.8\end{tabular}&\begin{tabular}[c]{@{}c@{}}2.5\end{tabular}&\begin{tabular}[c]{@{}c@{}}4.3\end{tabular}&\begin{tabular}[c]{@{}c@{}}5.9\end{tabular}&\begin{tabular}[c]{@{}c@{}}7.1\end{tabular}\\\cline{2-7}&SA/VIC&\begin{tabular}[c]{@{}c@{}}2.2\end{tabular}&\begin{tabular}[c]{@{}c@{}}0.8\end{tabular}&\begin{tabular}[c]{@{}c@{}}1.1\end{tabular}&\begin{tabular}[c]{@{}c@{}}1.5\end{tabular}&\begin{tabular}[c]{@{}c@{}}1.7\end{tabular}\\\hline
\end{tabular}
\end{table}

\subsubsection{Modal analysis} 
\begin{figure} 
\centering
\includegraphics[width=8.0cm, keepaspectratio]{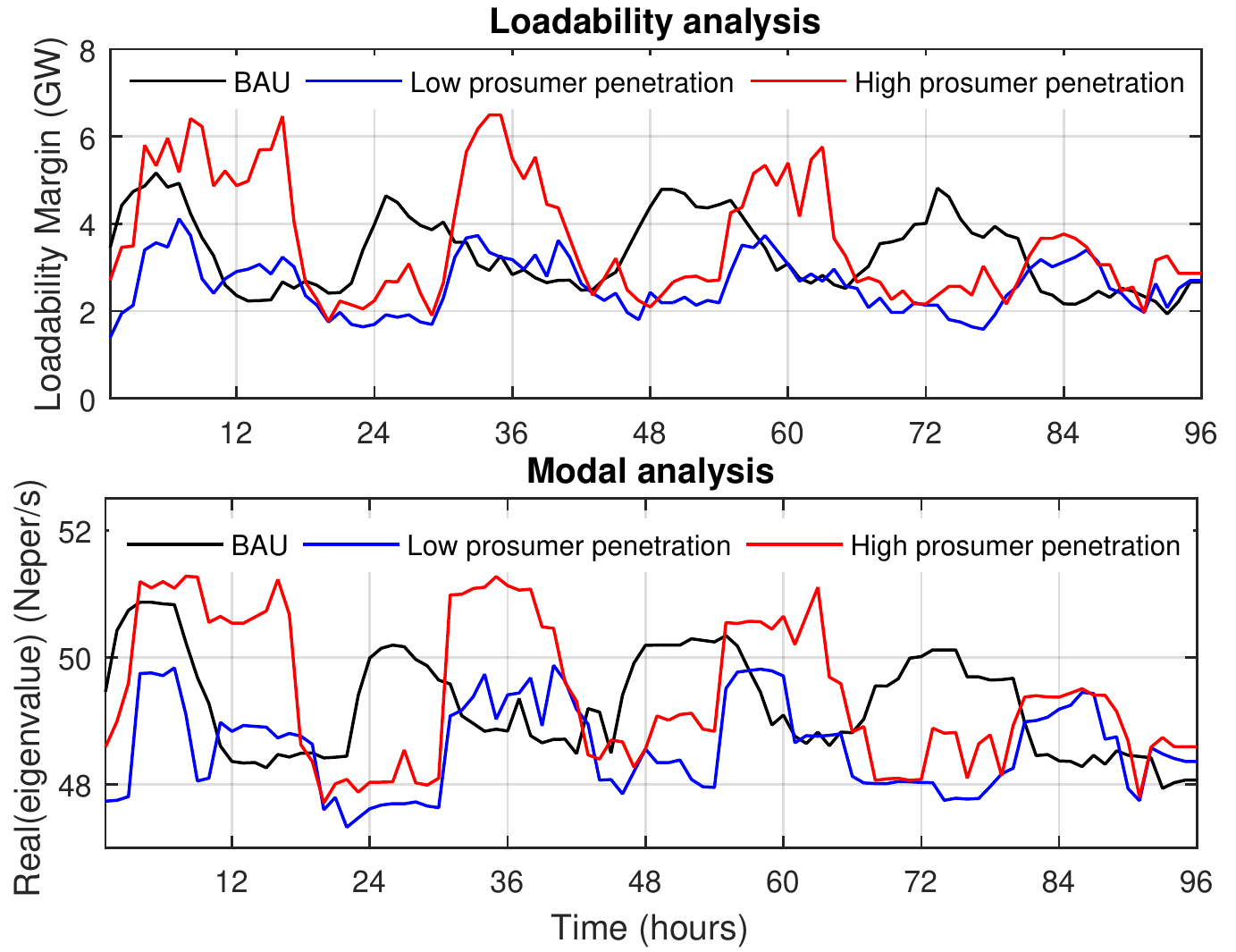}
\caption{Comparison of the loadability and modal analysis results for a typical summer week with high demand (12-15 January) for the load increase scenario NEM for BAU, and low and high prosumer penetration.}
\label{figure:LDBSA}
\end{figure}

Using the market dispatch results, modal analysis of the reduced V-Q sub-matrix of the power flow Jacobian is performed to assess voltage stability. The smallest real part of the V-Q sub-matrix' eigenvalues is used as a relative measure of the proximity to voltage instability. Furthermore, the associated eigenvectors provide information on the critical voltage modes and the weak points in the grid, that is, the areas that are most prone to voltage instability. The results are summarized in Table \ref{tab:eig results}.
With the increased RES penetration, the average of the minimum of the real part of all eigenvalues, henceforth called the minimum eigenvalue, is reduced from \SI{49.5}{Np \per \second} for the BAU scenario to \SI{42.3}{Np \per \second} for the ZP scenario. 
With an increased prosumer penetration, the average of the minimum eigenvalue over the simulated year increases from \SI{42.3}{Np \per \second} (ZP) to \SI{48.8}{Np \per \second} (HP). 
Observe in Fig. \ref{figure:LDBSA} that the results of the modal analysis confirm the results loadability analysis. The general trend remains the same; higher RES penetration with conventional demand reduces the minimum real value of an eigenvalue, which implies a lower voltage stability margin.

\begin{table}[]
\centering
\caption{Modal analysis results}
\label{tab:eig results}
\begin{tabular}{|c|c|c|c|c|c|}
\hline
& \multicolumn{5}{c|}{Scenarios}
\\\hline&BAU&ZP&LP&MP&HP\\\hline
\begin{tabular}[c]{@{}c@{}}Avg. of minimum\\ Real(eig) (Neper/s)\end{tabular}&\begin{tabular}[c]{@{}c@{}}49.5\end{tabular}&\begin{tabular}[c]{@{}c@{}}42.3\end{tabular}&\begin{tabular}[c]{@{}c@{}}44.8\end{tabular}&\begin{tabular}[c]{@{}c@{}}45.6\end{tabular}&\begin{tabular}[c]{@{}c@{}}48.4\end{tabular}\\\hline
\begin{tabular}[c]{@{}c@{}}Nodes with highest\\participation factor in\\critical voltage modes\end{tabular}&\begin{tabular}[c]{@{}c@{}}506\\306\\308\end{tabular}&\begin{tabular}[c]{@{}c@{}}506\\505\\410\end{tabular}&\begin{tabular}[c]{@{}c@{}}506\\505\\410\end{tabular}&\begin{tabular}[c]{@{}c@{}}505\\506\\410\end{tabular}&\begin{tabular}[c]{@{}c@{}}505\\410\\408\end{tabular}\\\hline
Unstable hours&0&175&37&12&0\\\hline
\end{tabular}
\end{table}

Another observation that can be made from the modal analysis concerns the location of the weak points in the system, that is, the buses with the highest participation factor in the critical voltage modes. These clearly change with the increased RES penetration.
In the BAU Scenario, the weakest points are typically buses with large loads (e.g. 306 and 308, representing Melbourne).
With the increased RES penetration and no prosumers (ZP), however, the weakest part of the system become buses located close to RESs (e.g. 505 and 410, representing, respectively a large wind farm in SA and a large solar PV farm in QLD).

Further, we observed that the voltage stability margin in the system improves significantly when there are more synchronous generators in the grid, due to their superior reactive power support capability compared to RES. The participation factor analysis of the renewable scenarios revealed that SA and QLD are the most voltage constrained regions where the penetration of WFs and utility PVs is higher compared to other regions in the NEM, which can be, to a large extent, mitigated with a sufficiently large penetration of prosumers. This clearly illustrates that RESs and prosumers change power system stability in ways that have not been experienced before, which requires a further in-depth analysis. 

\section{Conclusion}
The emergence of demand side technologies, in particular rooftop PV, battery storage and energy management systems is changing the way electricity consumers source and consume electric power, which requires new demand models for the long-term analysis of future grids. In this paper, we propose a generic demand model that captures the aggregate effect of a large number of prosumers on the load profile that can be used in market simulation.
The model uses a bi-level optimization framework, in which the upper level employs a unit commitment problem to minimize generation cost, and the lower level problem maximizes collective prosumers' self-consumption. To that effect, the model implicitly assumes an efficient mechanism for demand response aggregation, for example peer-to-peer energy trading or any other form of transactive energy. Given that any efficient mechanism for demand response aggregation that aims to minimize generation cost will utilize self-generation first, the self-consuming assumption appears reasonable at the level of abstraction assumed in long-term future grid scenario analysis. Moreover, the model is generic in that it does not depend on specific practical implementation details that will vary in the long-run.

To showcase the efficacy of the proposed model, we study the impact of prosumers on the performance, loadability and voltage stability of the Australian NEM with a high RES penetration. The results show that an increased prosumer penetration flattens the demand profile, which increases loadability and voltage stability, except in situations with a low underlying demand and an excess of RES generation, where the aggregate demand might increase due to battery charging. The analysis also revealed that with a high RES penetration, the weakest points in the network move from large load centers to areas with high RES penetration. Also, loadability and voltage stability are highly dependent on the amount of synchronous generation due to their superior reactive power capability compared to RES, which requires further analysis.

Our future work will focus on designing market mechanisms for large scale aggregation of distributed energy resources to enable the participation of prosumers in power system operation.

\appendix
\section*{MILP Formulation of the Lower-level KKT Optimality Conditions}
For the sake of completeness, we provide the complete optimization problem below.
The upper-level problem (\ref{eq:UC_utility})-(\ref{LC13}) remains the same.
The price-anticipatory prosumer behavior is captured in constraint (\ref{LC20}), defined as the solution to the lower-level prosumer aggregation problem. To be able to incorporate it in the upper-level problem, we first need to convert it into a set of MILP constraints. We do that by using the optimality (KKT) conditions of the underlying optimization problem. 
The KKT conditions follow from the associated Lagrangian:
\begin{align}
\mathcal{L}(\textbf{x}, \boldsymbol{\lambda}, \boldsymbol{\mu}) = f(\textbf{x})+ \boldsymbol{\lambda}^\top \textbf{g}(\textbf{x}) + \boldsymbol{\mu}^\top \textbf{h}(\textbf{x}),
\end{align}	
where $\textbf{x} = \{p_{\text{d}}^{\text{flx}}, p_{\text{b}}, e_{\text{b}} \}$ is the decision vector of the lower-level problem;  $f(\textbf{x})$ is the objective function (\ref{LC14}); $\textbf{g}(\textbf{x})$ is the vector of equality constraints (\ref{LLC2})-(\ref{LLC7}); $\textbf{h}(\textbf{x})$ is the vector of inequality constraints (\ref{LLC1})-(\ref{LLC6}); $\boldsymbol{\lambda}^\top = \{ \lambda_{h}^{m,p}, \lambda_{h}^{m,e} \}$ is the vector of Lagrange multipliers associated with the equality constraints, and $\boldsymbol{\mu}^\top = \{ \mu_{h}^{\text{flx},m}, \mu_{h}^{m,\underline{p}}, \mu_{h}^{m,\overline{p}}, \mu_{h}^{m,\underline{e}},\mu_{h}^{m,\overline{e}}\}$ is the vector of Lagrange multipliers associated with the inequality constraints.

The KKT optimality conditions consist of:
\begin{itemize}
\item primal feasibility (\ref{LLC2})-(\ref{LLC6}) with the associated Lagrange multipliers given in parentheses:
\end{itemize}
\begin{align}
& \qquad \qquad p_{\text{d},h}^{\text{flx},m} + p_{\text{pv},h}^{m} -  p_{\text{b},h}^{m} - p_{\text{d},h}^{\text{u},m}  =  0,   &( \lambda_{h}^{m,p} ) \label{LLC2}\\
& \qquad \qquad e_{\text{b},h}^{m} - \eta_\text{b}^m e_{\text{b},h-1}^{m} -  p_{\text{b},h}^{m}  =  0, & ( \lambda_{h}^{m,e} )   \label{LLC7}\\
& \qquad \qquad -p_{\text{d},h}^{\text{flx},m}  \leq 0,     &( \mu_{h}^{\text{flx},m} ) \label{LLC1}\\
& \qquad \qquad \underline{p}_{\text{b}}^{m} - p_{\text{b},h}^{m}  \leq 0,    &( \mu_{h}^{m,\underline{p}} )  \label{LLC3a}\\
& \qquad \qquad p_{\text{b},h}^{m} - \overline{p}_{\text{b}}^{m}  \leq 0,    &( \mu_{h}^{m,\overline{p}} )  \label{LLC4}\\
& \qquad \qquad \underline{e}_{\text{b}}^{m} -e_{\text{b},h}^{m}  \leq 0,    &( \mu_{h}^{m,\underline{e}} )  \label{LLC5}\\
& \qquad \qquad e_{\text{b},h}^{m} -\overline{e}_{\text{b}}^{m}  \leq 0,    &( \mu_{h}^{m,\overline{e}} )  \label{LLC6}
\end{align}
\begin{itemize}
\item dual feasibility:  
\end{itemize}
\begin{align}
& \mu_{h}^{\text{flx},m}, \mu_{h}^{m,\underline{p}}, \mu_{h}^{m,\overline{p}}, \mu_{h}^{m,\underline{e}}, \mu_{h}^{m,\overline{e}} \geq 0,
\end{align}	
\begin{itemize}
\item stationarity (\ref{KKT1})-(\ref{KKT3}):
\end{itemize}
\begin{align}
&\frac{\partial \mathcal{L}}{\partial p_{\text{d},h}^{\text{flx},m}} = 1 + \lambda_{h}^{m,p} - \mu_{h}^{\text{flx},m}  = 0, \label{KKT1}\\
&\frac{\partial \mathcal{L}}{\partial p_{\text{b},h}^{m}} = - \lambda_{h}^{m,p} - \lambda_{h}^{m,e} - \mu_{h}^{m,\underline{p}} + \mu_{h}^{m,\overline{p}}  = 0,\label{KKT2}\\
&\frac{\partial \mathcal{L}}{\partial e_{\text{b},h}^{m}} = \lambda_{h}^{m,e} -  \eta_\text{b}^m \lambda_{h+1}^{m,e} - \mu_{h}^{m,\underline{e}} + \mu_{h}^{m,\overline{e}}  = 0,\label{KKT3}
\end{align}	
\begin{itemize}
\item complementary slackness (\ref{SLC1})-(\ref{SLC5}):
\end{itemize}
\begin{align}
&p_{\text{d},h}^{\text{flx},m} \mu_{h}^{\text{flx},m} = 0, \label{SLC1}\\
&(-\underline{p}_{\text{b}}^{m} + p_{\text{b},h}^{m}) \mu_{h}^{m,\underline{p}} = 0,  \label{SLC2}\\
&(-p_{\text{b},h}^{m} + \overline{p}_{\text{b}}^{m}) \mu_{h}^{m,\overline{p}} = 0, \label{SLC3}\\
&(-\underline{e}_{\text{b}}^{m} + e_{\text{b},h}^{m}) \mu_{h}^{m,\underline{e}} = 0, \label{SLC4}\\
&(-e_{\text{b},h}^{m} + \overline{e}_{\text{b}}^{m}) \mu_{h}^{m,\overline{e}} = 0. \label{SLC5}
\end{align}

Complementary slackness conditions~\eqref{SLC1}-\eqref{SLC5} are bilinear, however they can be linearized by introducing one binary variable $b$ and a sufficiently large and positive constant $M$ resulting in two constraints per complementarity slackness condition (\ref{Eq_middle})-(\ref{Eq_last}).

Given that the KKT conditions are sufficient and necessary for optimality, the bi-level constraint~\eqref{LC20} can be replaced with mixed integer linear constraints~\eqref{Eq_fir}-\eqref{Eq_last}:
\begin{align}
& p_{\text{d},h}^{\text{flx},m} + p_{\text{pv},h}^{m} -  p_{\text{b},h}^{m} - p_{\text{d},h}^{\text{u},m}  =  0, \label{Eq_fir}\\ 
& e_{\text{b},h}^{m} - \eta_\text{b}^m e_{\text{b},h-1}^{m} -  p_{\text{b},h}^{m}  =  0 \label{eLLC7},\\
&  1 - \mu_{h}^{\text{flx},m} + \lambda_{h}^{m,p}    = 0,\label{eKKT1}\\
&  - \lambda_{h}^{m,p} - \mu_{h}^{m,\underline{p}} + \mu_{h}^{m,\overline{p}} - \lambda_{h}^{m,e}   = 0,\label{eKKT2}\\
&  - \mu_{h}^{m,\underline{e}} + \mu_{h}^{m,\overline{e}} + \lambda_{h}^{m,e} -  \eta_\text{b}^m \lambda_{h-1}^{m,e}   = 0,\label{eKKT3}\\
& \underline{p}_{\text{b}}^{m} - p_{\text{b},h}^{m}   \leq 0,  \label{eLLC3}\\
& p_{\text{b},h}^{m} -\overline{p}_{\text{b}}^{m}  \leq 0, \label{eLLC4}\\
& \underline{e}_{\text{b}}^{m} - e_{\text{b},h}^{m}   \leq 0  , \label{eLLC5}\\
& e_{\text{b},h}^{m} -\overline{e}_{\text{b}}^{m}  \leq 0  , \label{eLLC6}\\
& p_{\text{d},h}^{\text{flx},m} - M{^\text{flx}} b_{h}^{\text{flx},m}  \leq 0, \label{Eq_middle}\\
& \mu_{h}^{\text{flx},m} - M^\text{flx} (1-b_{h}^{\text{flx},m})  \leq 0,\\
& p_{\text{b},h}^{m} - M^{\underline{p}} b_{h}^{m,\underline{p}}  \leq 0,\\
& \mu_{h}^{m,\underline{p}} - M^{\underline{p}} (1-b_{h}^{m,\underline{p}})  \leq 0,\\
& -p_{\text{b},h}^{m} -M^{\overline{p}} b_{h}^{m,\overline{p}}   \leq 0,\\
& \mu_{h}^{m,\overline{p}} - M^{\overline{p}} (1-b_{h}^{m,\overline{p}})   \leq 0 ,\\
& e_{\text{b},h}^{m} - M^{\underline{e}} b_{h}^{m,\underline{e}}  \leq 0,\\
& \mu_{h}^{m,\underline{e}} - M^{\underline{e}} (1-b_{h}^{m,\underline{e}})   \leq 0,\\
& -e_{\text{b},h}^{m} - M^{\overline{e}} b_{h}^{m,\overline{e}}  \leq 0,\\
& \mu_{h}^{m,\overline{e}}  -M^{\overline{e}} (1-b_{h}^{m,\overline{e}})  \leq 0\label{Eq_last},
\end{align} 
which adds the following additional decision variables to the upper-level problem (\ref{eq:UC_utility})-(\ref{LC13}):
\begin{align*}
& b_{h}^{\text{flx},m},b_{h}^{m,\underline{p}},b_{h}^{m,\overline{p}},b_{h}^{m,\underline{e}},b_{h}^{m,\overline{e}} \in \{0,1\},\\
& p_{\text{b},h}^{m}, \lambda_{h}^{m,p}, \lambda_{h}^{m,e}  \in \mathbb{R}, \\
& p_{\text{d},h}^{\text{flx},m}, e_{\text{b},h}^{m}, \mu_{h}^{\text{flx},m}, \mu_{h}^{m,\underline{p}}, \mu_{h}^{m,\overline{p}}, \mu_{h}^{m,\underline{e}},\mu_{h}^{m,\overline{e}} \in \mathbb{R}_+.
\end{align*}
The resulting optimization problem is an MILP that can be solved efficiently with off-the-shelf solvers.
\ifCLASSOPTIONcaptionsoff
  \newpage
\fi

\bibliographystyle{IEEEtran}
\bibliography{GenDemModel}

\begin{IEEEbiography}[{\includegraphics[width=1in,height=1.25in,clip,keepaspectratio]{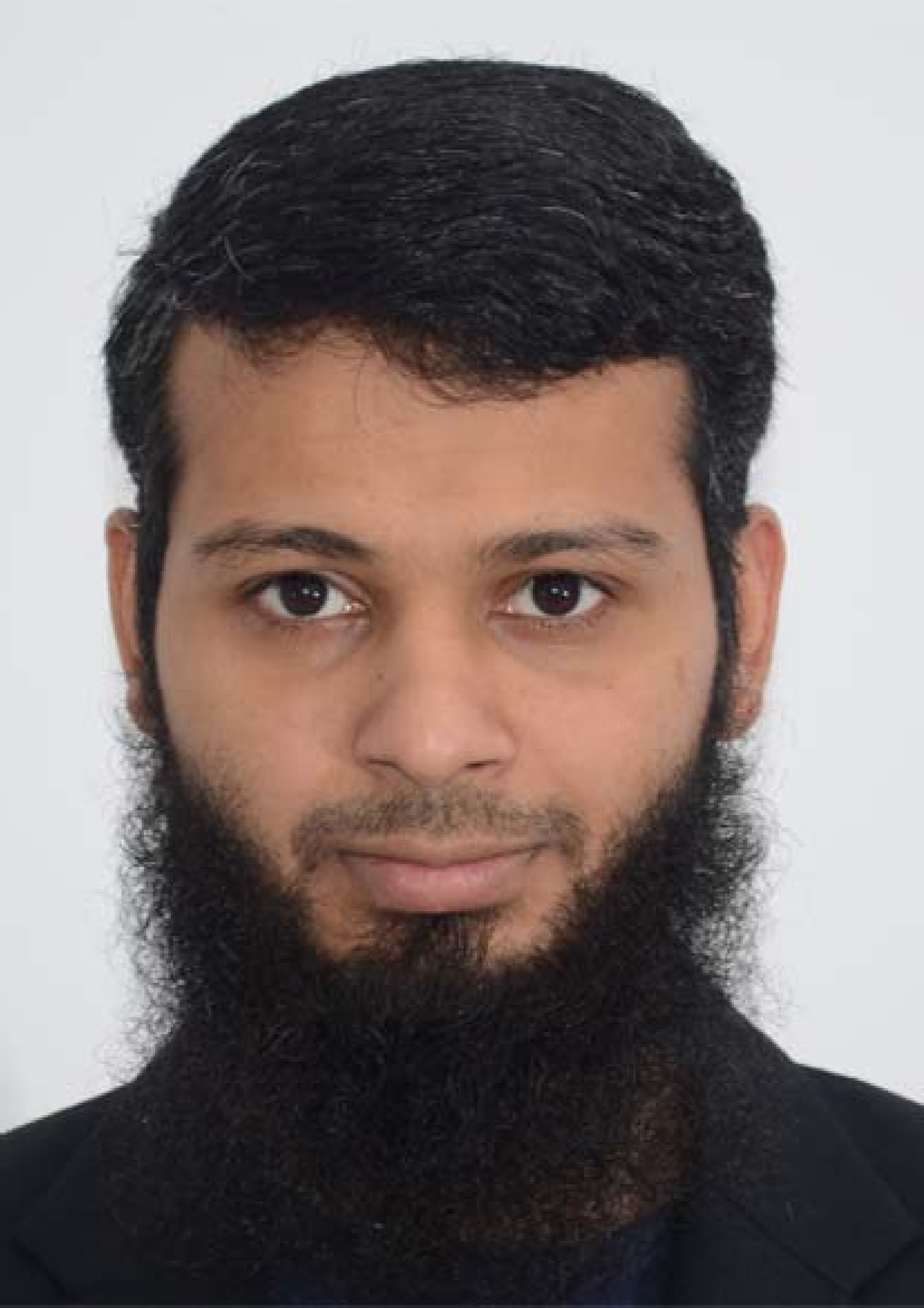}}]{Shariq Riaz} (S'14) received the B.Sc. and M.Sc. degrees in electrical engineering from the University of Engineering and Technology Lahore, Lahore, Pakistan in 2009 and 2012, respectively. He is currently pursuing the Ph.D. degree from The University of Sydney, Sydney, Australia. His expertise is in power system operation, demand response, and electricity markets. His current research interests include integration of renewable generation, storage and prosumers into electricity markets, demand response and operation of concentrated solar thermal plants.
\end{IEEEbiography}

\begin{IEEEbiography}[{\includegraphics[width=1in,height=1.25in,clip,keepaspectratio]{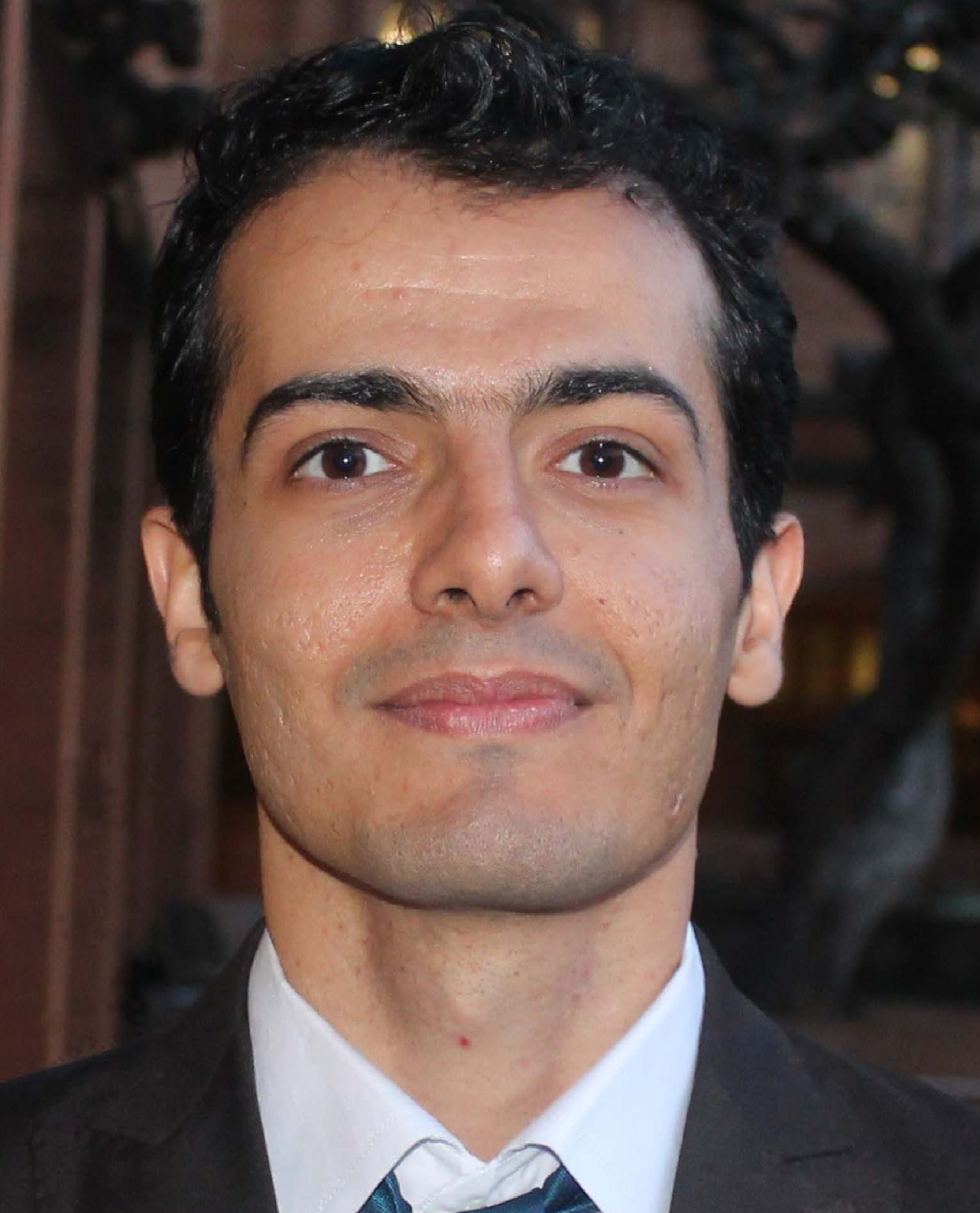}}]{Hesamoddin Marzooghi} (S'13-M'16) graduated with Bachelor and Masters degrees in electrical engineering (power) both with first class honors from Shiraz University, Iran in 2009 and 2011, respectively. He was awarded his PhD degree in electrical engineering (power) at The University of Sydney in 2016. He is now a research associate at The Central Queensland University (CQUniversity), Australia, and an academic visitor at The University of Manchester, UK. Hesamoddin's research interests are in stability analysis and control of power systems with high penetration of renewable energy sources, distributed generation and energy storage, planning of future grids, and applications of intelligent systems in power engineering.
\end{IEEEbiography}

\begin{IEEEbiography}[{\includegraphics[width=1in,height=1.25in,clip,keepaspectratio]{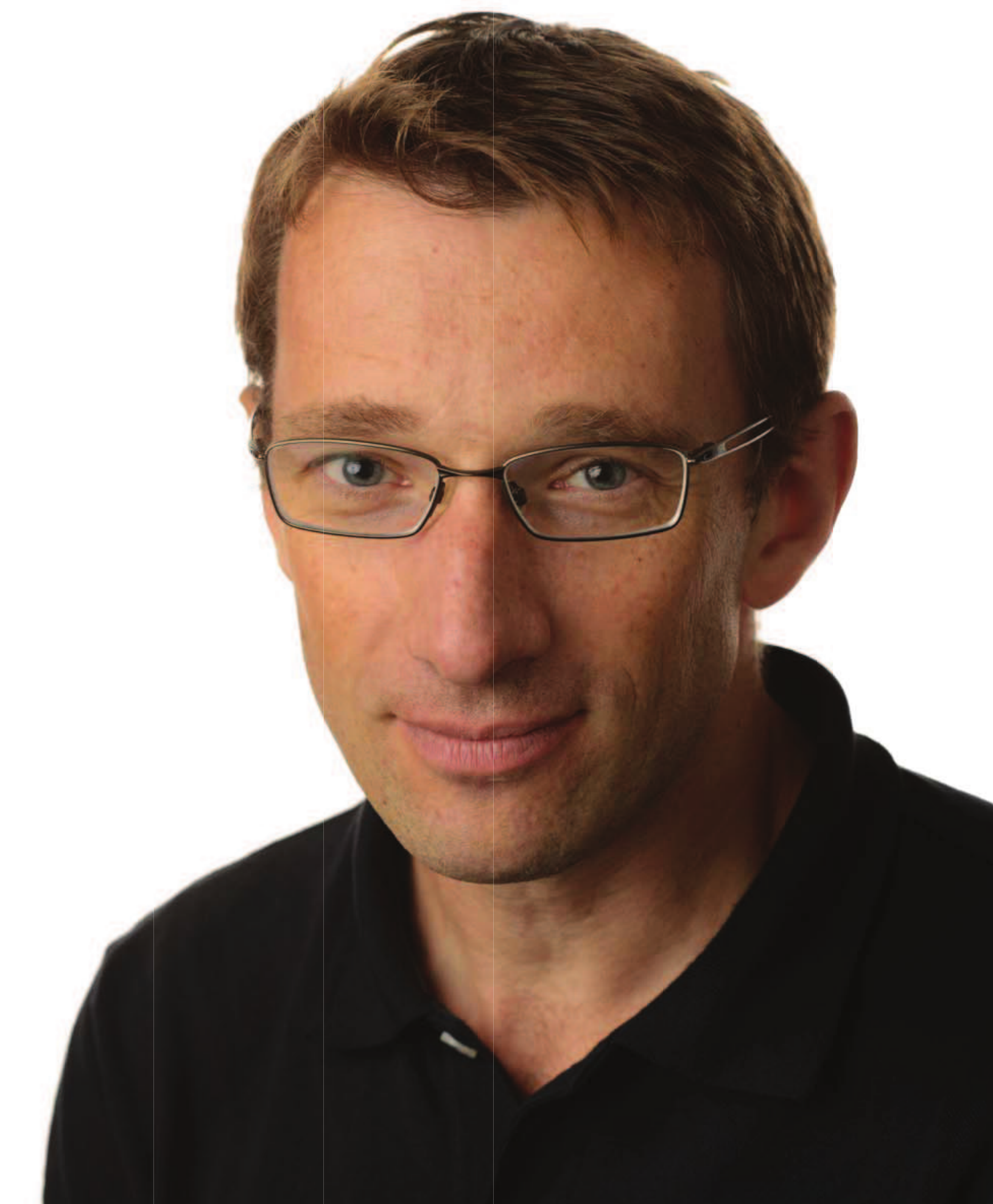}}]{Gregor Verbi\v{c}} (S'98-M’03-SM'10) received the B.Sc., M.Sc., and Ph.D. degrees in electrical engineering from the University of Ljubljana, Ljubljana, Slovenia, in 1995, 2000, and 2003, respectively. In 2005, he was a NATO-NSERC Postdoctoral Fellow with the University of Waterloo, Waterloo, ON, Canada. Since 2010, he has been with the School of Electrical and Information Engineering, The University of Sydney, Sydney, NSW, Australia. His expertise is in power system operation, stability and control, and electricity markets. His current research interests include integration of renewable energies into power systems and markets, optimization and control of distributed energy resources, demand response, and energy management in residential buildings. He was a recipient of the IEEE Power and Energy Society Prize Paper Award in 2006. He is an Associate Editor of the IEEE Transactions on Smart Grid.
\end{IEEEbiography}

\begin{IEEEbiography}[{\includegraphics[width=1in,height=1.25in,clip,keepaspectratio]{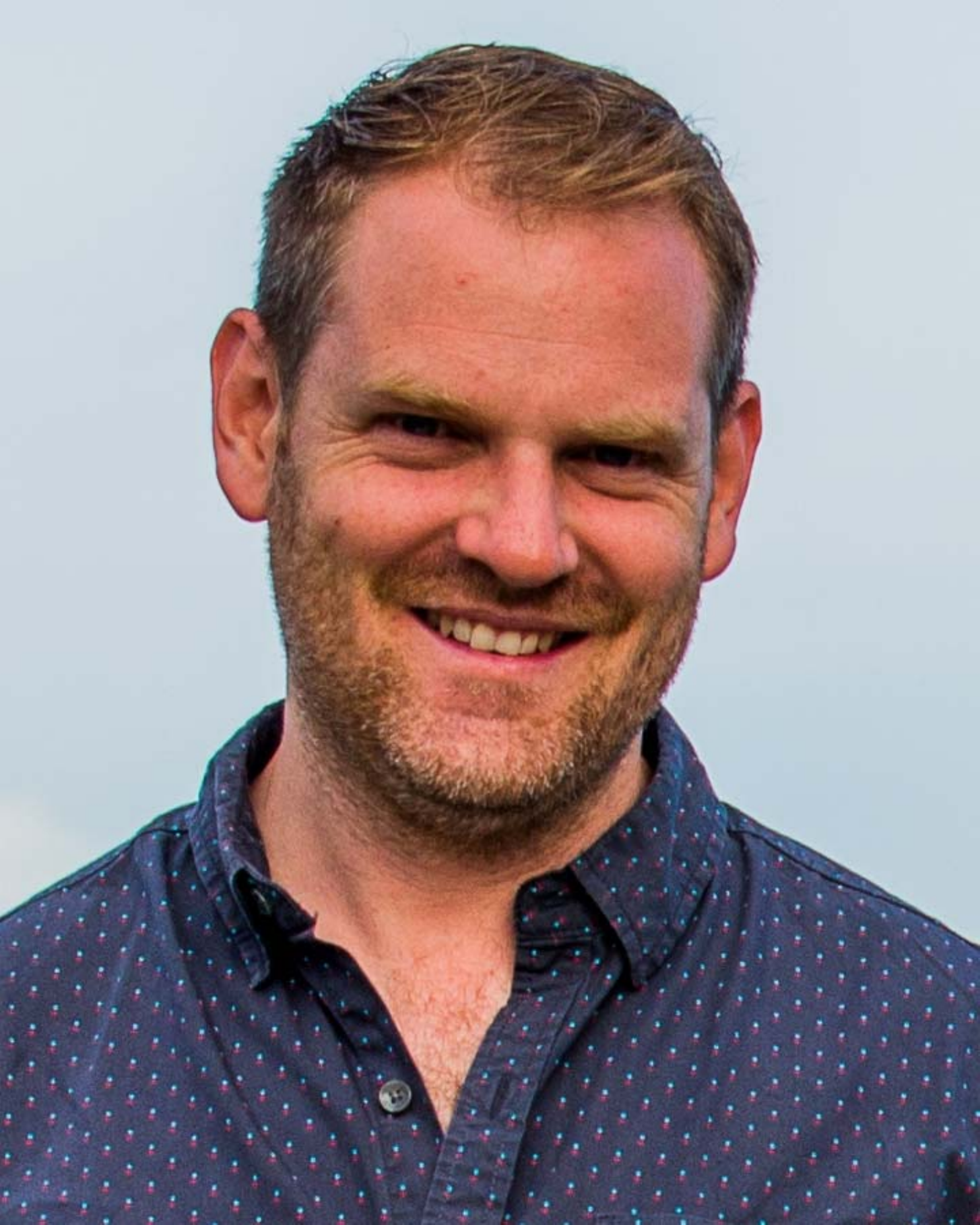}}]{Archie C. Chapman} (M'14) received the B.A. degree in math and political science, and the B.Econ. (Hons.) degree from the University of Queensland, Brisbane, QLD, Australia, in 2003 and 2004, respectively, and the Ph.D. degree in computer science from the University of Southampton, Southampton, U.K., in 2009. He is currently a Research Fellow in Smart Grids with the School of Electrical and Information Engineering, Centre for Future Energy Networks, The University of Sydney, Sydney, NSW, Australia. He has expertise is in game-theoretic and reinforcement learning techniques for optimization and control in large distributed systems. His research focuses on integrating renewables into legacy power networks, using distributed energy and load scheduling methods, and on designing tariffs and market mechanisms that support efficient use of existing infrastructure and new controllable devices. 
\end{IEEEbiography}

\begin{IEEEbiography}[{\includegraphics[width=1in,height=1.25in,clip,keepaspectratio]{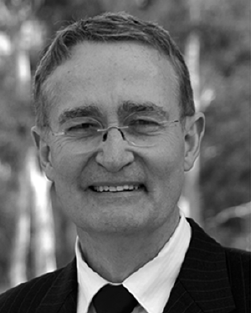}}]{David J. Hill} (S'72-M'76-SM'91-F'93-LF'14) received the Ph.D. degree in electrical engineering from the University of Newcastle, Australia, in 1976. He holds the Chair of Electrical Engineering in the Department of Electrical and Electronic Engineering at the University of Hong Kong. He is also a part-time Professor of Electrical Engineering at The University of Sydney, Australia. During 2005–2010, he was an Australian Research Council Federation Fellow at the Australian National University. Since 1994, he has held various positions at the University of Sydney, Australia, including the Chair of Electrical Engineering until 2002 and again during 2010–2013 along with an ARC Professorial Fellowship. He has also held academic and substantial visiting positions at the universities of Melbourne, California (Berkeley), Newcastle (Australia), Lund (Sweden), Munich, Hong Kong (City and Polytechnic). His general research interests are in control systems, complex networks, power systems and stability analysis. His work is now mainly on control and planning of future energy networks and basic stability and control questions for dynamic networks. 

Prof. Hill is a Fellow of the Society for Industrial and Applied Mathematics, the Australian Academy of Science, and the Australian Academy of Technological Sciences and Engineering. He is also a Foreign Member of the Royal Swedish Academy of Engineering Sciences.
\end{IEEEbiography}

\end{document}